\documentclass[10pt,a4paper]{article}
\usepackage[cp1251]{inputenc}
\usepackage{amsfonts, amssymb, amsthm, amsmath}
\usepackage{mathtext}
\usepackage{longtable}
\usepackage[russian]{babel}
\usepackage{color}

\usepackage[left=1cm,right=1cm,top=2cm,bottom=2cm]{geometry}
\usepackage{graphicx}

\theoremstyle{remark}

\renewcommand{\r}{\mathbb R}

\renewcommand{\le}{\leqslant}
\renewcommand{\ge}{\geqslant}
\newcommand{\I}{\mathbb{I}}

\renewcommand{\phi}{\varphi}

\newcommand{\eqd}{\stackrel{d}{=}}
\newcommand{\tod}{\stackrel{d}{\longrightarrow}}

\renewcommand{\refname}{References}
\title{A note on functional limit
theorems for compound Cox processes\thanks{Research supported by
Russian Scientific Foundation, project 14-11-00364.}}
\author{V. Yu. Korolev\thanks{Faculty of
Computational Mathematics and Cybernetics, Lomonosov Moscow State
University; Institute of Informatics Problems FRC CSC RAS;
victoryukorolev@yandex.ru}, A. V. Chertok\thanks{Faculty of
Computational Mathematics and Cybernetics, Lomonosov Moscow State
University, Euphoria Group LLC; a.v.chertok@gmail.com}, A. Yu.
Korchagin\thanks{Faculty of Computational Mathematics and
Cybernetics, Lomonosov Moscow State University;
sasha.korchagin@gmail.com}, E. V. Kossova\thanks{Higher School of
Economics National Research University, Moscow;
e$\_$kossova@mail.ru}, A. I. Zeifman\thanks{Vologda State
University; Institute of Informatics Problems FRC CSC RAS;  ISEDT
RAS; a$\_$zeifman@mail.ru} }

\date{}

\begin{document}

\maketitle

{\small

{\bf Abstract:} An improved version of the functional limit theorem
is proved establishing weak convergence of random walks generated by
compound doubly stochastic Poisson processes (compound Cox
processes) to L{\'e}vy processes in the Skorokhod space under more
realistic moment conditions. As corollaries, theorems are proved on
convergence of random walks with jumps having finite variances to
L{\'e}vy processes with variance-mean mixed normal distributions, in
particular, to stable L{\'e}vy processes, generalized hyperbolic and
generalized variance-gamma L{\'e}vy processes.

\smallskip


}

\section{Introduction}

In financial mathematics the evolution of (the logarithms of) stock
prices and financial indexes on small time horizons is often modeled
by random walks. The simplest example of such an approach is the
Cox--Ross--Rubinstein model (see, e. g., \cite{Shiryaev1999}). At
the same time most successful (adequate) models of the dynamics of
(the logarithms of) financial indexes on large time horizons are
subordinated Wiener processes (processes of Brownian motion with
random time such as generalized hyperbolic processes, in particular,
variance gamma (VG) processes and
normal$\backslash\!\backslash$inverse Gaussian (NIG) processes.
These models very well describe the observed heavy-tailedness and
leptokurticity of the empirical (statistical) distributions of the
increments of financial indices and, in particular, of stock prices
on comparatively short time intervals.

Functional limit theorems are a quite natural link between random
walks and subordinated Wiener processes. The operation of
subordination gives a good explanation of the presence of heavy
tails in the empirical distributions of the increments of (the
logarithms of) stock prices and financial indexes.

In the book \cite{GnedenkoKorolev1996} and the papers
\cite{Korolev1997, Korolev2000} it was proposed to model the
evolution of non-homogeneous chaotic stochastic processes, in
particular, of the dynamics of stock prices and financial indexes,
by random walks generated by compound doubly stochastic Poisson
processes (compound Cox pocesses). A {\it doubly stochastic Poisson
process} (also called a {\it Cox process}) is a stochastic point
process of the form $N_1(\Lambda(t))$, where $N_1(t)$, $t\geq0$, is
a homogeneous Poisson process with unit intensity and the stochastic
process $\Lambda(t)$, $t\geq0$, is independent of $N_1(t)$ and
possesses the following properties: $\Lambda(0)=0$, ${\sf
P}(\Lambda(t)<\infty)=1$ for any $t>0$, the sample paths of
$\Lambda(t)$ do not decrease and are right-continuous. A compound
Cox process is a random sum of independent identically distributed
random variables in which the number of summands follows a Cox
process. Similar continuous-time random walks were considered in
\cite{gz97, Huang2012}.

This approach, based on the universal principle of non-decrease of
entropy in closed systems, was developed in \cite{BeningKorolev2002,
KorolevBeningShorgin2011, Korolev2011}. In the paper
\cite{KorolevChertokKorchaginZeifman2015} this approach was
successfully applied to modeling the evolution of limit order books
in the high-frequency financial trading systems. In the framework of
this approach the principal idea is that the moments at which the
state of the system under consideration changes form a {\it chaotic}
point stochastic process on the time axis. Moreover, this point
process turns out to be non-stationary (time-non-homogeneous)
because the changes of the state of the limit order book are to a
great extent subject to the influence of non-stationary information
flows. As is known, most reasonable probabilistic models of
non-stationary (time-non-homogeneous) chaotic point processes are
{\it doubly stochastic Poisson processes} also called {\it Cox
processes} (see, e. g., \cite{Grandell1976, BeningKorolev2002}).
These processes are defined as Poisson processes with stochastic
intensities. Pure Poisson processes can be regarded as best models
of stationary (time-homogeneous) chaotic flows of events
\cite{BeningKorolev2002}. Recall that the attractiveness of a
Poisson process as a model of homogeneous discrete stochastic chaos
is due to at least two reasons. First, Poisson processes are point
processes characterized by that time intervals between successive
points are independent random variables with one and the same
exponential distribution and, as is well known, the exponential
distribution possesses the maximum differential entropy among all
absolutely continuous distributions concentrated on the nonnegative
half-line with finite expectations, whereas the entropy is a natural
and convenient measure of uncertainty. Second, the points forming
the Poisson process are uniformly distributed along the time axis in
the sense that for any finite time interval $[t_1,t_2]$, $t_1<t_2$,
the conditional joint distribution of the points of the Poisson
process which fall into the interval $[t_1,t_2]$ under the condition
that the number of such points is fixed and equals, say, $n$,
coincides with the joint distribution of the order statistics
constructed from an independent sample of size $n$ from the uniform
distribution on $[t_1,t_2]$ whereas the uniform distribution
possesses the maximum differential entropy among all absolutely
continuous distributions concentrated on finite intervals and very
well corresponds to the conventional impression of an absolutely
unpredictable random variable (see, e. g.,
\cite{GnedenkoKorolev1996, BeningKorolev2002}).

In \cite{KorolevZaksZeifman2013a} some functional limit theorems
were proved establishing convergence of random walks generated by
compound Cox processes with jumps possessing finite variances to
L{\'e}vy processes with symmetric distributions including symmetric
strictly stable L{\'e}vy processes. In the paper
\cite{KorolevChertokKorchaginZeifman2015} these results were
extended to a non-symmetric case and applied to modeling the
evolution of the order flow imbalance process, an integral
characteristic of the behavior of the limit order book.

The present paper presents a further development of the models and
techniques proposed in our previous papers
\cite{KorolevZaksZeifman2013a, KorolevChertokKorchaginZeifman2015}.
In this paper we improve and generalize the results of the mentioned
papers by relaxing the conditions and correcting some inaccuracies.
Our presentation essentially relies on the techniques developed in
\cite{KorolevZaksZeifman2013a}.

The paper is organized as follows. Section 2 contains some basic
material on the Skorokhod space, stable distributions and L{\'e}vy
processes. In section 3 we prove general functional limit theorem
establishing the conditions for convergence of compound Cox
processes to L{\'e}vy processes in the Skorokhod space in terms of
the behavior of the cumulative intensities of Cox processes. For
this purpose we extend the classical results presented, say, in
\cite{JacodShiryaev2003}. In Section 4 we consider the conditions
for the convergence of compound Cox processes with elementary jumps
possessing finite variances to the L{\'e}vy processes with
variance-mean mixed normal one-dimensional distributions, that is,
to subordinated Wiener processes.

\section{Skorokhod space. L{\'e}vy processes}

Let $D=D[0,1]$ be the space of real functions defined on $[0,1]$,
right-continuous and having finite left-side limits.

Let $\mathcal{F}$ be the class of strictly increasing mappings of
$[0,1]$ onto itself. Let $f$ be a non-decreasing function on $[0,1]$
with $f(0)=0$, $f(1)=1$. Set $\|f\|=\sup_{s\neq
t}\left|\log\left[\big(f(t)-f(s)\big)/(t-s)\right]\right|$. If
$\|f\|<\infty$, then the function $f$ is continuous and strictly
increasing and, hence, belongs to $\mathcal{F}$.

Define the distance $d_0(x,y)$ in the set $D[0,1]$ as the greatest
lower bound of the set of positive numbers $\epsilon$, for which
$\mathcal{F}$ contains a function $f$ such that $\|f\|\le\epsilon$
and $\sup_t|x(t)-y(f(t))|\le\epsilon$.

It can be shown that the space $D[0,1]$ is complete with respect to
the distance $d_0$. The metric space $\mathcal{D}=(D[0,1],d_0)$ is
called {\it the Skorokhod space}. Everywhere in what follows we will
consider stochastic processes as $\mathcal{D}$-valued random
elements.

Let $X, X_1,X_2,...$ be $\mathcal{D}$-valued random elements. Let
$T_X$ be a subset of $[0,1]$ such that $0\in T_X$, $1\in T_X$ and if
$0<t<1$, then $t\in T_X$ if and only if ${\sf P}\bigl(X(t)\neq
X(t-)\bigr)=0$. The following theorem establishing sufficient
conditions for the weak convergence of stochastic processes in
$\mathcal{D}$ (denoted below as $\Longrightarrow$ and assumed as
$n\to\infty$) is well-known.

\smallskip

{\sc Theorem A.} {\it Let}
$\bigl(X_n(t_1),...,X_n(t_k)\bigr)\Longrightarrow
\bigl(X(t_1),...,X(t_k)\bigr)$ {\it for any natural $k$ and
$t_1,...,t_k$ belonging to $T_X$. Let ${\sf P}\bigl(X(1)\neq
X(1-)\bigr)=0$ and let there exist a non-decreasing continuous
function $F$ on $[0,1]$, such that for any $\epsilon>0$}
$$
{\sf P}\bigl(|X_n(t)-X_n(t_1)|\ge\epsilon,\
|X_n(t_2)-X_n(t)|\ge\epsilon\bigr)
\le\epsilon^{-2\nu}\bigl[F(t_2)-F(t_1)\bigr]^{2\gamma}\eqno(1)
$$
{\it for $t_1\le t\le t_2$ and $n\ge1$, where $\nu\ge 0$,
$\gamma>1/2$. Then} $X_n\Longrightarrow X$.

\smallskip

The {\sc proof} of Theorem A can be found, for example, in
\cite{Billingsley1968}.

\smallskip

Everywhere in what follows the symbol $\eqd$ stands for the
coincidence of distributions.

By a L{\'e}vy process we will mean a stochastic process $X(t)$,
$t\ge0$, possessing the following properties: (i) $X(0)=0$ almost
surely; (ii) $X(t)$ is the process with independent increments, that
is, for any $N\ge1$ and $t_0,t_1,...,t_N$ ($0\le t_0\le t_1\le...\le
t_N$) the random variables $X(t_0)$, $X(t_1)-X(t_0)$, $...,$
$X(t_N)-X(t_{N-1})$ are jointly independent; (iii) $X(t)$ is a
homogeneous process, that is, $X(t+h)-X(t)\eqd X(s+h)-X(s)$ for any
$s,t,h>0$; (iv) the process $X(t)$ is stochastically continuous,
that is, for any $t\ge0$ and $\epsilon>0$ $ \lim_{s\to t}{\sf
P}(|X(t)-X(s)|>\epsilon)=0$; (v) sample paths of the process $X(t)$
are right-continuous and have finite left-side limits.

Denote the characteristic function of the random variable $X(t)$ as
$\psi_t(s)$ ($\psi_t(s)={\sf E}e^{isX(t)}$, $s\in\mathbb{R}$). The
following statement describes a well-known property of L{\'e}vy
processes.

\smallskip

{\sc Lemma 1.} {\it Let $X=X(t)$, $t\ge 0$, be a L{\'e}vy process.
For any $t>0$ the characteristic function of the random variable
$X(t)$ is infinitely divisible and has the form
$$
\psi_{t}(s) = \big[\psi_{1}(s)\big]^{t} = \big[{\sf
E}\,e^{isX(1)}\big]^t,\ \ \ \  s\in\mathbb{R}.\eqno(2)
$$
Conversely, let $Y$ be an arbitrary infinitely divisible random
variable. Then the family of infinitely divisible distributions with
characteristic functions of the form $\big[{\sf E}\,e^{isY}\big]^t$
completely determines finite-dimensional distributions of a L{\'e}vy
process $X(t)$, $t\ge 0$, moreover,} $X(1)\stackrel{d}{=}Y$.

\smallskip

The properties of L{\'e}vy processes are described in detail in
\cite{Bertoin1996, Sato1999}.

By $G_{\alpha,\theta}(x)$ we will denote the distribution function
of the strictly stable law with the characteristic exponent $\alpha$
and parameter $\theta$ corresponding to the characteristic function
$\mathfrak{g}_{\alpha,\theta}(s)=\exp\big\{-|s|^{\alpha}\exp\big\{-\frac{i}{2}
\pi\theta\alpha\mathrm{sign}s\big\}\big\}$, $s\in\r$, where
$0<\alpha\le2$,
$|\theta|\le\theta_{\alpha}=\min\{1,\frac{2}{\alpha}-1\}$. To
symmetric strictly stable distributions there corresponds the value
$\theta=0$. To one-sided stable distributions there correspond the
values $\theta=1$ and $0<\alpha\le1$.

In what follows any random variable with the distribution function
$G_{\alpha,\theta}(x)$, $0<\alpha<2$, will be denoted
$Z_{\alpha,\theta}$. It is well known that ${\sf
E}|Z_{\alpha,\theta}|^{\delta}<\infty$ for any
$\delta\in(0,\alpha)$, but the moments of orders higher or equal to
$\alpha$ of the random variable $Z_{\alpha,\theta}$ do not exist
(see, e. g., \cite{Zolotarev1983}).

The distribution function of the standard normal law ($\alpha=2$,
$\theta=0$) will be denoted $\Phi(x)$,
$\Phi(x)=\int_{-\infty}^x\phi(z)dz,$
$\phi(x)=\frac{1}{\sqrt{2\pi}}e^{-x^2/2}$.

It is well known that the distribution function $G_{\alpha,0}(x)$ of
the symmetric strictly stable law can be represented as a scale
mixture of normal laws:
$$
G_{\alpha,0}(x)=\int_{0}^{\infty}\Phi\Big(\frac{x}{\sqrt{u}}\Big)dG_{\alpha/2,1}(u),\
\ \ x\in\mathbb{R}
$$
(see, e.g., \cite{Zolotarev1983}, Theorem 3.3.1). To this
representation there corresponds the following relation in terms of
characteristic functions:
$$
\mathfrak{g}_{\alpha,0}(s)=\int_{0}^{\infty}\exp\Big\{-\frac{s^2u}{2}\Big\}dG_{\alpha/2,1}(u),\
\ \ s\in\mathbb{R}.%\eqno(1)
$$

A L{\'e}vy process $X(t)$, $t\ge0$, will be called {\it
$\alpha$-stable}, if ${\sf P}\big(X(1)<x\big)=G_{\alpha,\theta}(x)$,
$x\in\mathbb{R}$. It can be shown (see, e.g.,
\cite{EmbrechtsMaejima2002}) that if $X(t)$, $t\ge0$, is a L{\'e}vy
process, then $X(t)$ is $\alpha$-stable if and only if
$$
X(t)\eqd t^{1/\alpha}X(1),\ \ \ t\ge0.\eqno(3)%\eqno(5)
$$

\section{Convergence of compound Cox processes to L{\'e}vy processes}

In what follows without noticeable loss of generality we will
consider stochastic processes defined for $0\le t\le1$. Actually,
this means that we consider the behavior of compound Cox processes
on finite time horizons. The equality of the right bound of the
horizon to one can be achieved by an appropriate choice of the units
of measurement of time. In other words, we will concentrate on
studying the case of the Skorokhod space $\mathcal{D}$.

In order to introduce reasonable asymptotics which formalizes the
condition of <<infinite>> growth of intensities of the flows of
informative events, and makes it possible to construct asymptotic
(<<heavy-traffic>>) approximations to Cox processes, fix a time
instant $t$ and introduce an auxiliary parameter $n$. Everywhere in
what follows the convergence will be meant as $n\to\infty$ unless
otherwise specified. So, consider a sequence of compound Cox
processes of the form
$$
Q_n(t)=\sum\nolimits_{i=1}^{N^{(n)}_{1}(\Lambda_n(t))}X_{n,i},\ \ \
t\ge0,\eqno(4)
$$
where $\{N^{(n)}_{1}(t),\, t\ge0\}_{n\ge1}$ is a sequence of Poisson
processes with unit intensities; for each $n=1,2,...$ the random
variables $X_{n,1},X_{n,2},...$ are identically distributed; for any
$n\ge1$ the random variables $X_{n,1},X_{n,2},...$ and the process
$N^{(n)}_{1}(t)$, $t\ge0$, are independent; for each $n=1,2,...$
$\Lambda_n(t)$, $t\ge0$, is a subordinator, that is, a
non-decreasing positive L{\'e}vy process, independent of the process
$$
Z_n(t)=\sum\nolimits_{i=1}^{N^{(n)}_{1}(t)}X_{n,i},\ \ \ t\ge
0,\eqno(5)
$$
and such that $\Lambda_n(0)=0$.

Assume that there exist $\delta\in(0,1]$, $\delta_1\ge\frac12$ and
the constants $C_n\in(0,\infty)$ providing for all $t\in(0,1]$ the
validity of the inequality
$$
{\sf E}\Lambda^{\delta}_n(t)\le (C_nt)^{\delta_1}.\eqno(6)
$$

For example, assume that $\Lambda_n(t)$ is a stable L{\'e}vy
process, that is, ${\sf P}\big(\Lambda_n(1)<x\big)=G_{\alpha,1}(x)$
with some $0<\alpha<1$. Then for any
$\rho\in(0,\alpha)\subseteq(0,1]$ we have ${\sf
E}\Lambda_n^{\rho}(1)<\infty$ and, moreover, in accordance with (3)
we have ${\sf E}\Lambda_n^{\rho}(t)=t^{\rho/\alpha}{\sf
E}\Lambda_n^{\rho}(1)$, that is, condition (6) holds for any
$\delta\in[\alpha/2,\alpha)$ with
$\delta_1=\delta/\alpha\in[\frac12,1)$ and $C_n=\big({\sf
E}\Lambda_n^{\delta}(1)\big)^{1/\alpha}$.

Assume that
$$
0<m_n^{\beta}\equiv{\sf E}|X_{n,1}|^{\beta}<\infty.\eqno(7)
$$
for some $\beta\in(0,1]$.

Everywhere in what follows for definiteness we assume that
$\sum\nolimits_{i=1}^0=0$.

From (4) and (5) it is easy to see that $Q_n(t)=Z_n(\Lambda_n(t))$.
Since for each $n\ge1$ both $Z_n(t)$ and $\Lambda_n(t)$ are
independent L{\'e}vy processes, and, moreover, $\Lambda_n(t)$ is a
subordinator, then the superposition $Q_n(t)=Z_n(\Lambda_n(t))$ is
also a L{\'e}vy process (see, e. g., Theorem 3.1.1 in
\cite{Kashcheev2001}). Hence the following statement follows.

\smallskip

{\sc Lemma 2}. {\it For any $0\le t_1<t_2<\infty$ and any $n\ge1$ we
have} $Q_n(t_2)-Q_n(t_1)\eqd Q_n(t_2-t_1)$.

\smallskip

{\sc Lemma 3}. {\it Let $Q_n(t)$ be a compound Cox process $(4)$
satisfying conditions $(6)$ and $(7)$. Then for any $t\in[0,1]$ and
any $\epsilon>0$ we have} ${\sf P}\big(|Q_n(t)|\ge\epsilon\big)\le
(\epsilon^{-\beta}m_n^{\beta})^{\delta}\cdot\big(C_nt\big)^{\delta_1}$.

\smallskip

{\sc Proof}. Since one-dimensional distributions of the Cox process
(4) are mixed Poisson, we have
$$
{\sf P}\big(|Q_n(t)|\ge\epsilon\big)= {\sf
P}\bigg(\Big|\sum\nolimits_{j=1}^{N^{(n)}_1(\Lambda_n(t))}X_{n,j}\Big|\ge\epsilon\bigg)=
\sum\nolimits_{k=0}^{\infty}{\sf
P}\big(N^{(n)}_1(\Lambda_n(t))=k\big){\sf
P}\Big(\Big|\sum\nolimits_{j=1}^kX_{n,j}\Big|\ge\epsilon\Big)=
$$
%$$
%=\sum\nolimits_{k=0}^{\infty}\Big(\int_{0}^{\infty}e^{-\lambda}\frac{\lambda^k}{k!}\,d{\sf
%P}\big(\Lambda_n(t)<\lambda\big)\Big){\sf
%P}\Big(\Big|\sum\nolimits_{j=1}^kX_{n,j}\Big|\ge\epsilon\Big)=
%$$
$$
=\int_{0}^{\infty}\Big[\sum\nolimits_{k=0}^{\infty}e^{-\lambda}\frac{\lambda^k}{k!}{\sf
P}\Big(\Big|\sum\nolimits_{j=1}^kX_{n,j}\Big|\ge\epsilon\Big)\Big]\,d{\sf
P}\big(\Lambda_n(t)<\lambda\big).\eqno(8)
$$
The change of the order of summation and integration is possible due
to the obvious uniform convergence of the series. Continue (8) by
the sequential application of the Markov and Jensen inequalities
with $\delta\in(0,1]$ taking part in (6) and $\beta\in(0,1]$ taking
part in (7). As a result we obtain
$$
{\sf P}\big(|Q_n(t)|\ge\epsilon\big)\le
\frac{1}{\epsilon^{\beta\delta}}\int_{0}^{\infty}\Big[\sum\nolimits_{k=0}^{\infty}e^{-\lambda}\frac{\lambda^k}{k!}{\sf
E}\Big|\sum\nolimits_{j=1}^kX_{n,j}\Big|^{\beta\delta}\Big]\,d{\sf
P}\big(\Lambda_n(t)<\lambda\big) \le
$$
$$
\le
\frac{1}{\epsilon^{\beta\delta}}\int_{0}^{\infty}\Big[\sum\nolimits_{k=0}^{\infty}e^{-\lambda}\frac{\lambda^k}{k!}\Big({\sf
E}\Big|\sum\nolimits_{j=1}^kX_{n,j}\Big|^{\beta}{\Big)\!}^{\delta}\,\Big]\,d{\sf
P}\big(\Lambda_n(t)<\lambda\big),\eqno(9)
$$
since with $\delta\in(0,1]$ the function $f(x)=x^{\delta}$ is
concave for $x\ge0$. It is easy to see that ${\sf
E}\big|\sum\nolimits_{j=1}^kX_{n,j}\big|^{\beta}\le\sum\nolimits_{j=1}^k{\sf
E}|X_{n,j}|^{\beta}=km_n^{\beta}$ for $0<\beta\le 1$. Therefore,
continuing (9) with the account of the Jensen inequality for concave
functions and (6), we obtain
$$
{\sf P}\big(|Q_n(t)|\ge\epsilon\big)\le
\frac{(m_n^{\beta})^{\delta}}{\epsilon^{\beta\delta}}\int_{0}^{\infty}\!\!\Big(\sum\nolimits_{k=0}^{\infty}e^{-\lambda}\frac{k^{\delta}\lambda^k}{k!}\Big)d{\sf
P}\big(\Lambda_n(t)<\lambda\big)=
\frac{(m_n^{\beta})^{\delta}}{\epsilon^{\beta\delta}}\int_{0}^{\infty}\!\!{\sf
E}\big[N^{(1)}_1(\lambda)\big]^{\delta}\,d{\sf
P}\big(\Lambda_n(t)<\lambda\big)\le
$$
$$
\le
\frac{(m_n^{\beta})^{\delta}}{\epsilon^{\beta\delta}}\int_{0}^{\infty}\!\!\big[{\sf
E}N^{(1)}_1(\lambda)\big]^{\delta}\,d{\sf
P}\big(\Lambda_n(t)<\lambda\big) =
\frac{(m_n^{\beta})^{\delta}}{\epsilon^{\beta\delta}}\int_{0}^{\infty}\!\!\lambda^{\delta}\,d{\sf
P}\big(\Lambda_n(t)<\lambda\big)=
\frac{(m_n^{\beta})^{\delta}}{\epsilon^{\beta\delta}}\cdot{\sf
E}\Lambda_n^{\delta}(t)\le
\frac{(m_n^{\beta})^{\delta}}{\epsilon^{\beta\delta}}\cdot\big(C_nt\big)^{\delta_1}.
$$
The lemma is proved.

\smallskip

To establish weak convergence of the stochastic processes $Q_n(t)$
in the Skorokhod space $\mathcal{D}$, first it is required to find
the limit distribution of the random variables $Q_n(t)$ for each
$t>0$. The symbol $\tod$ will denote convergence in distribution,
that is, pointwise convergence of the distribution functions in all
continuity points of the limit distribution function.

Let $t=1$. Denote $N_n=N^{(n)}_{1}(\Lambda_n(1))$. Assume that for
some $k_n\in\mathbb{N}$ the convergence
$$
{\sf P}(X_{n,1}+...+X_{n,k_n}<x)\tod H(x)\eqno(10)
$$
takes place, where $H(x)$ is some infinitely divisible distribution
function.

Also assume that
$$
{\sf P}\big(\Lambda_n(1)<k_n x\big)\tod {\sf P}(U<x),\eqno(11) %(13)
$$
where $U$ is a nonnegative random variable such that its
distribution is not degenerate in zero. Notice that since
$\Lambda_n(t)$ is a L{\'e}vy process, then the random variable $U$
is infinitely divisible being the weak limit of infinitely divisible
random variables.

\smallskip

{\sc Lemma 4}. {\it Let $N_n=N_1^{(n)}(\Lambda_n)$, $n\ge1$, where
$\{N_1^{(n)(t)},\ t\ge0\}$, $n=1,2,\ldots$ are standard Poisson
processes and $\Lambda_n$, $n=1,2,\ldots$ are positive random
variables such that for each $n\ge1$ the random variable $\Lambda_n$
is independent of the process $N_1^{(n)}(t)$. Then ${\sf
P}(N_n<k_nx)\tod A(x)$ for some infinitely increasing sequence $k_n$
of real numbers and some distribution function $A(x)$ if and only
if} ${\sf P}(\Lambda_n<k_nx)\tod A(x)$.

\smallskip

For the proof see \cite{GnedenkoKorolev1996}.

\smallskip

From Lemma 4 it follows that convergence (11) is equivalent to
$$
{\sf P}\big(N_n<k_n x\big)\tod {\sf P}(U<x).\eqno(12)
$$
By the Gnedenko--Fahim transfer theorem \cite{GnedenkoFahim1969}
conditions (10) and (12) imply that
$$
Q_n(1)=X_{n,1}+...+X_{n,N_n}\tod Q,\eqno(13)
$$
where $Q$ is a random variable with the characteristic function
$\mathfrak{f}(s)=\int_{0}^{\infty}\big(h(s)\big)^u\,d{\sf P}(U<u)$,
$h(s)$ being the characteristic function corresponding to the
distribution function $H(x)$. Note that the distribution function
$H(x)$  may not satisfy the condition $H(-x)=1-H(x)$ for all
$x\ge0$, that is, it may not be symmetric.

Let $Y$ be an infinitely divisible random variable with the
distribution function $H(x)$. Since both $Y$ and $U$ are infinitely
divisible, we can define independent L{\'e}vy processes $Y(t)$ and
$U(t)$, $t\ge0$, such that $Y(1)\eqd Y$ and $U(1)\eqd Y$. Then with
the account of Lemma 1 it is easy to verify that
$\mathfrak{f}(s)={\sf E}e^{isQ}={\sf E}\exp\big\{isY(U(1))\big\}$,
$s\in\mathbb{R}$, that is, $Q\eqd Y(U(1))$. Moreover, repeating the
reasoning from \cite{Kashcheev2001} (see Theorem 3.3.1 there), we
can easily see that the random variable $Q$ is infinitely divisible
and hence, we can define a L{\'e}vy process $Q(t)$, $t\ge0$, such
that $Q(1)\eqd Q$. From Lemma 1 and the abovesaid it follows that we
can regard $Q(t)$ as the superposition: $Q(t)\eqd Y(U(t))$.

Since according to (13) we have
$Q_n(1)=\sum\nolimits_{i=1}^{N_n}X_{n,i}\Longrightarrow Q(1)$, and
both $Q_n(t)$ and $Q(t)$ are L{\'e}vy processes, then, using (2) we
can conclude that for any $t>0$
$$
Q_n(t)=\sum\nolimits_{i=1}^{N_{n,1}(\Lambda_n(t))}X_{n,i}\tod Q(t).
\ \ \ \eqno(14)
$$

Since the processes $Q_n(t)$ and $Q(t)$, $0\le t\le 1$, are L{\'e}vy
processes, then almost all their sample paths belong to the
Skorokhod space ${\cal D}$.

Consider the question what additional conditions are required to
provide the weak convergence of the compound Cox process $Q_n(t)$ to
the L{\'e}vy process $Q(t)$ in the space ${\cal D}$. We will
consider each of the conditions of Theorem A one by one.

First, without loss of generality, let $0\le t_1<t_2<...<t_k \le 1$.
The convergence $\bigl(Q_n(t_1),...,Q_n(t_k)\bigr)\tod
\bigl(Q(t_1),...,Q(t_k)\bigr)$ is equivalent to the convergence
$$\bigl(Q_n(t_1),Q_n(t_2)-Q_n(t_1),...,Q_n(t_k)-Q_n(t_{k-1})\bigr)
\tod
$$
$$
\tod\bigl(Q(t_1),Q(t_2)-Q(t_1),...,Q(t_k)-Q(t_{k-1})\bigr),\eqno(15) %(22)
$$
since the linear transform $(x_1,x_2,...,x_{k-1},x_k)\longmapsto
(x_1,x_2-x_1,...,x_k-x_{k-1})$ of $\mathbb{R}^k$ to $\mathbb{R}^k$
is one-to-one and continuous in both directions. But convergence
(15) follows from (14) and the fact that both $Q_n(t)$ and $Q(t)$
are L{\'e}vy processes.

Second, we have to check the condition ${\sf P}\bigl(Q(1)\neq
Q(1-)\bigr)=0$. This condition holds if and only if
$\lim_{t\to1-}{\sf P}\bigl(|Q(1)-Q(t)|>\epsilon\bigr)=0$ for any
$\epsilon>0$ (see relation (15.16) in \cite{Billingsley1968}).
Consider ${\sf P}\bigl(|Q(1)-Q(t)|>\epsilon\bigr)$. Since $Q(t)$ is
a L{\'e}vy process, then $Q(1)-Q(t)\eqd Q(1-t)$ by Lemma 2.
Therefore, ${\sf P}\bigl(|Q(1)-Q(t)|>\epsilon\bigr)={\sf
P}\bigl(|Q(1-t)|>\epsilon\bigr)$. For each $\epsilon>0$ and each
$t\in[0,1]$ there exists an $\epsilon_t\in[\epsilon/2,\epsilon]$
such that the points $\pm\epsilon_t$ are continuity points of the
distribution function of the random variable $Q(1-t)$. Since
$Q_n(t)\tod Q(t)$ for each $t\in[0,1]$, then $ {\sf
P}\bigl(|Q(1-t)|>\epsilon_t\bigr)=\lim_{n\to\infty}{\sf
P}\bigl(|Q_n(1-t)|>\epsilon_t\bigr)$. Thus, for any $\epsilon>0$ and
any $t\in[0,1]$ we have
$$
{\sf P}\bigl(|Q(1-t)|>\epsilon\bigr)\le{\sf
P}\bigl(|Q(1-t)|>\epsilon_t\bigr)=\lim_{n\to\infty}{\sf
P}\bigl(|Q_n(1-t)|>\epsilon_t\bigr)\le\limsup_{n\to\infty}{\sf
P}\bigl(|Q_n(1-t)|>\epsilon_t\bigr).\eqno(16) %(23)
$$
Continuing (16) with the account of (6) and applying Lemma 3, for
$\delta\in(0,1]$ taking part in (6) we obtain
$$
{\sf P}\bigl(|Q(1-t)|>\epsilon\bigr)\le\limsup_{n\to\infty}{\sf
P}\bigl(|Q_n(1-t)|>\epsilon_t\bigr)\le
$$
$$
\le
\limsup_{n\to\infty}\big(\epsilon_t^{-\beta}m_n^{\beta}\big)^{\delta}\big(C_n|1-t|\big)^{\delta_1}\le
\big((2/\epsilon)^{\beta\delta}|1-t|\big)^{\delta_1}\limsup_{n\to\infty}\,(m_n^{\beta})^{\delta}C_n^{\delta_1}.\eqno(17) %(24)
$$
Therefore, if
$$
K\equiv\limsup_{n\to\infty}C_n^{\delta_1/\delta}m_n^{\beta}<\infty,\eqno(18) %(25)
$$
then (17) implies
$$
\lim_{t\to1-}{\sf P}\bigl(|Q(1)-Q(t)|>\epsilon\bigr)\le
4(K\epsilon^{-\beta})^{\delta}\lim_{t\to1-}|1-t|^{\delta_1}=0.
$$
Third, check condition (1) under the assumption that (6) and (18)
hold. As it has been noted above, $Q_n(t)$ is a L{\'e}vy process and
hence, it has independent increments. Therefore,
$$
{\sf
P}\big(|Q_n(t)-Q_n(t_1)|\ge\epsilon,\,|Q_n(t_2)-Q_n(t)|\ge\epsilon\big)=
 {\sf P}\big(|Q_n(t)-Q_n(t_1)|\ge\epsilon\big)\cdot{\sf
P}\big(|Q_n(t_2)-Q_n(t)|\ge\epsilon\big).\eqno(19) %(26)
$$
Consider the first multiplier on the right-hand side of (19). By
Lemma 2, $Q_n(t)-Q_n(t_1)\eqd Q_n(t-t_1)$. With the account of (18),
by Lemma 3 we obtain
$$
{\sf P}\big(|Q_n(t)-Q_n(t_1)|\ge\epsilon\big)={\sf
P}\big(|Q_n(t-t_1)|\ge\epsilon\big)\le(K\epsilon^{-\beta})^{\delta}|t-t_1|^{\delta_1}.\eqno(20) %(27)
$$
For the second multiplier on the right-hand side of (19) we
similarly obtain
$$
{\sf P}\big(|Q_n(t_2)-Q_n(t)|\ge\epsilon\big)={\sf
P}\big(|Q_n(t_2-t)|\ge\epsilon\big)\le(K\epsilon^{-\beta})^{\delta}|t_2-t|^{\delta_1}.\eqno(21) %(28)
$$
Thus, from (20) and (21) it follows that
$$
{\sf
P}\big(|Q_n(t)-Q_n(t_1)|\ge\epsilon,\,|Q_n(t_2)-Q_n(t)|\ge\epsilon\big)\le
(K\epsilon^{-\beta})^{2\delta}\big[(t-t_1)(t_2-t)\big]^{\delta_1}\eqno(22) %(29)
$$
It is easy to see that for any $t_1\le t\le t_2$ we have
$(t-t_1)(t_2-t)\le{\textstyle\frac14}(t_2-t_1)^2$. Substituting this
estimate in (22) we obtain ${\sf
P}\big(|Q_n(t)-Q_n(t_1)|\ge\epsilon,\,|Q_n(t_2)-Q_n(t)|\ge\epsilon\big)\le
\epsilon^{-2\beta\delta}\big[\frac12{K(t_2-t_1)}\big]^{2\delta_1}$.
Therefore, if conditions (6) and (18) hold, then condition (1) holds
with $F(t)\equiv \frac12Kt$, $\nu=\beta\delta$ and
$\gamma=\delta_1$.

Summarizing this reasoning related to checking the conditions of
Theorem A, we arrive at the following statement.

\smallskip

{\sc Theorem 1.} {\it Assume that the random variables
$\{X_{n,j}\}_{j\ge1}$, $n=1,2,...$, $($the jumps of the compound Cox
process $Q_n(t)$, see $(4))$, satisfy conditions $(10)$ with some
$k_n\in\mathbb{N}$ and $(7)$ with some $\beta\in(0,1]$. Let the
processes $Q_n(t)$ be lead by non-decreasing positive L{\'e}vy
processes $\Lambda_n(t)$ satisfying conditions $(6)$ with some
$\delta\in(0,1]$, $\delta_1\ge\frac12$ and $(11)$ with the same
$k_n$. Also assume that condition $(18)$ holds. Then the processes
$Q_n(t)$ weakly converge in the Skorokhod space $\mathcal{D}$ to the
L{\'e}vy process $Q(t)$ such that
$$
{\sf E}\exp\{isQ(1)\}=\int_{0}^{\infty}\big(h(s)\big)^u\,d{\sf
P}(U<u),\ \ \ s\in\mathbb{R},\eqno(23) %(30)
$$
where $h(s)$ is the characteristic function corresponding to the
distribution function $H(x)$ in $(10)$.}

\smallskip

It is worth noting that actually Theorem 1 deals with the
well-studied weak convergence of special semimartingales with
stationary increments, see, e. g., \cite{JacodShiryaev2003}.
However, the superposition-type structure of the processes
considered in the present paper makes it possible to relax the
conditions required in the general case, say, in Corollary VII.3.6
of \cite{JacodShiryaev2003} where it is assumed that (in our
terminology) $\delta=\delta_1=1$. Moreover, in
\cite{KorolevZaksZeifman2013a} and some subsequent papers a more
restrictive condition was used instead of (18).

From Theorem 1 it obviously follows that if in (10)
$H(x)=G_{\alpha,\theta}(x)$ for some admissible $\alpha\in(0,2)$ and
$\theta\in[-1,1]$ and the processes $\Lambda_n(t)$ are
asymptotically degenerate, that is, in (11) ${\sf P}(U=1)=1$, then
in (23) ${\sf E}\exp\{isQ(1)\}=\mathfrak{g}_{\alpha,\theta}(s)$,
that is, the limiting process is the stable L{\'e}vy process.

However, in \cite{KorolevZaksZeifman2013a} it was demonstrated that
stable L{\'e}vy processes can appear as limits for compound Cox
processes even when the variances of elementary increments of a
compound Cox process are finite. To generalize this result, consider
the following corollary of Theorem 1.

\smallskip

{\sc Corollary 1.} {\it Assume that the random variables
$\{X_{n,j}\}_{j\ge1}$, $n=1,2,...$, satisfy condition $(10)$ with
some $k_n\in\mathbb{N}$ and $H(x)=G_{\alpha,0}(x)$ for some
$\alpha\in(0,2]$ so that condition $(7)$ holds for any
$\beta\in(0,\alpha)\bigcap \,(0,1]$. Let the compound Cox processes
$Q_n(t)$ be lead by non-decreasing positive L{\'e}vy processes
$\Lambda_n(t)$ satisfying conditions $(11)$ with the same $k_n$ and
${\sf P}(U<x)=G_{\alpha',1}(x)$ for some $\alpha'\in(0,1]$ so that
condition $(6)$ holds with some $\delta\in[\alpha'/2,\alpha')$ and
$\delta_1=\delta/\alpha'\in[\frac12,1)$. Also assume that condition
$(18)$ holds for some $\beta\in(0,\alpha)$. Then the processes
$Q_n(t)$ weakly converge in the Skorokhod space $\mathcal{D}$ to the
L{\'e}vy process $Q(t)$ such that} ${\sf
P}\big(Q(1)<x\big)=G_{\alpha\alpha',0}(x)$, $x\in\r$.
\smallskip
\indent {\sc Proof.} To prove this result is suffices to notice that
in the case under consideration
$\mathfrak{f}_{\alpha,0}(s)=e^{-|s|^{\alpha}}$, $s\in\r$, so that in
(23)
$$
{\sf
E}\exp\{isQ(1)\}=\int_{0}^{\infty}e^{-u|s|^{\alpha}}dG_{\alpha',1}(u)=
%$$
%$$
%=
\int_{0}^{\infty}e^{-|su^{1/\alpha}|^{\alpha}}dG_{\alpha',1}(u)=
{\sf E}\exp\big\{isZ_{\alpha,0}Z_{\alpha',1}^{1/\alpha}\big\},\ \ \
s\in\r,
$$
where the random variables $Z_{\alpha,0}$ and $Z_{\alpha',1}$ are
independent. But from Theorem 3.3.1 in \cite{Zolotarev1983} it
follows that $Z_{\alpha,0}Z_{\alpha',1}^{1/\alpha}\eqd
Z_{\alpha\alpha',0}$. The corollary is proved.

\smallskip

In turn, the result from \cite{KorolevZaksZeifman2013a} mentioned
above follows from Corollary 1, if $\alpha=2$. This case corresponds
to the situation in which the variances of the summands (elementary
jumps) are assumed finite. As it has been already said, in most
applied problems there are no reasons to reject this assumption.
Therefore in what follows we will concentrate our attention on the
case where the elementary increments of the compound Cox process
have finite variances and consider the conditions of convergence of
compound Cox processes to some popular models, in particular, to
L{\'e}vy processes with variance-mean mixed normal one-dimensional
distributions such as generalized hyperbolic L{\'e}vy processes or
generalized variance-gamma L{\'e}vy processes.

\section{L{\'e}vy processes with
variance-mean mixed normal distributions as asymptotic
approximations to compound Cox processes}

Denote $a_n={\sf E}X_{n,1}$ and $\sigma_n^2={\sf D}X_{n,1}$. From
the classical theory of limit theorems it is well known that if, as
$n\to\infty$, the conditions
$$
k_na_n\longrightarrow a,\ \ k_n\sigma_n^2\longrightarrow \sigma^2 \
\mbox{ and } \ k_n{\sf
E}(X_{n,1}-a_n)^2\I(|X_{n,1}-a_n|\ge\epsilon)\longrightarrow
0\eqno(24) %(12)
$$
hold for some $a\in\mathbb{R}$, $0<\sigma^2<\infty$ and any
$\epsilon>0$, then convergence (10) takes place with
$H(x)\equiv\Phi\big(\sigma^{-1}(x-a)\big)$. In this case the
distribution function $F(x)$ of the limit random variable $Q(1)$ in
Theorem 1 is a variance-mean mixture of normal laws. Recently it was
demonstrated that normal variance-mean mixtures appear as limiting
in simple limit theorems for random sums of independent identically
distributed random variables \cite{Korolev2013}. Namely, let
$\{\xi_{n,j}\}_{j\ge1},$ $n=1,2,\ldots,$ be a double array of
row-wise (for each fixed $n$) independent and identically
distributed random variables. Let $\{\nu_n\}_{n\ge1}$ be a sequence
of integer nonnegative random variables such that for each $n\ge1$
the random variables $\nu_n,\xi_{n,1},\xi_{n,2},\ldots$ are
independent. Denote $ S_{n,k}=\xi_{n,1}+\ldots +\xi_{n,k}$. The
following theorem was proved in \cite{Korolev2013}.

\smallskip

{\sc Theorem B}. {\it Assume that there exist$:$ a sequence
$\{k_n\}_{n\ge1}$ of natural numbers and finite numbers
$\alpha\in\mathbb{R}$ and $\sigma>0$ such that
$$
{\sf P}\big(S_{n,k_n}<x\big)\tod
\Phi\Big(\frac{x-\alpha}{\sigma}\Big).\eqno(25)
$$
Assume that $\nu_n\to\infty$ in probability. Then the distribution
functions of random sums $S_{n,\nu_n}$ converge to some distribution
function $F(x):$ $ {\sf P}\big(S_{n,\nu_n}<x\big)\tod F(x)$, if and
only if there exists a distribution function $A(x)$ such that
$$
A(0)=0,\ \ \ F(x)=\int_{0}^{\infty}\Phi\Big(\frac{x-\alpha
z}{\sigma\sqrt{z}}\Big)dA(z),
$$
and ${\sf P}(\nu_n<xk_n)\tod A(x)$. }

\smallskip

Theorem B and Lemma 3 yield the following result.

\smallskip

{\sc Theorem 2.} {\it Assume that the random variables
$\{X_{n,j}\}_{j\ge1}$, $n=1,2,...$, $($the jumps of the compound Cox
process $Q_n(t)$, see $(4))$ possess finite variances and satisfy
conditions $(24)$ with some $k_n\in\mathbb{N}$, $a\in\mathbb{R}$ and
$\sigma>0$. Let the the processes $Q_n(t)$ be lead by non-decreasing
positive L{\'e}vy processes $\Lambda_n(t)$ satisfying condition
$(6)$ with some $\delta\in(0,1]$, $\delta_1\ge\frac12$. Also assume
that
$$
K\equiv\limsup_{n\to\infty}C_n^{\delta_1/\delta}(\sigma_n+|a_n|)<\infty.\eqno(26)
$$
Then the processes $Q_n(t)$ weakly converge in the Skorokhod space
$\mathcal{D}$ to a L{\'e}vy process $Q(t)$ if and only if there
exists a nonnegative random variable $U$ such that
$$
{\sf
P}\big(Q(1)<x\big)=\int_{0}^{\infty}\Phi\Big(\frac{x-au}{\sigma\sqrt{u}}\Big)d{\sf
P}(U<u),\ \ \ x\in\mathbb{R},\eqno(27) %(24)
$$
and condition $(11)$ holds with the same $k_n$.}

\smallskip

{\sc Proof.} It suffices to take $\beta=1$ in the proof of lemma 3
and notice that $m_n^1\le\sigma_n+|a_n|$ so that condition (26) can
play the role of (18).

\smallskip

The class of distributions of form (27) was systematically
considered by O. Barndorff-Nielsen and his colleagues \cite{BN1978,
BN1979, BN1982} in order to introduce {\it generalized hyperbolic
distributions} and study their properties.

The class of normal variance-mean mixtures (27) is very wide. For
example, it contains generalized hyperbolic laws with generalized
inverse Gaussian mixing distributions \cite{BN1977, BN1978},
generalized variance gamma (GVG) distributions with generalized
gamma mixing distributions \cite{KorolevSokolov2012,
KorolevZaks2013}, symmetric strictly stable laws with strictly
stable mixing distributions concentrated on the positive half-line,
generalized exponential power distributions, and many other types.

Generalized hyperbolic distributions demonstrate exceptionally high
adequacy when they are used to describe statistical regularities in
the behavior of characteristics of various complex open systems, in
particular, turbulent systems and financial markets. There are
dozens of dozens of publications dealing with models based on
generalized hyperbolic distributions. Recently it was discovered
that generalized variance gamma distributions demonstrate even
better fit to empirical data. Therefore below we will concentrate
our attention on functional limit theorems establishing the
convergence of compound Cox processes to generalized hyperbolic
L{\'e}vy processes and generalized variance gamma L{\'e}vy processes
yielding the possibility of the use of such processes as convenient
<<heavy-traffic>> {\it asymptotic} approximations.

Denote the density of the {\it generalized inverse Gaussian
distribution} by $p_{GIG}(x;\nu,\mu,\lambda)$,
$$
p_{GIG}(x;\nu,\mu,\lambda)=\frac{\lambda^{\nu/2}}{2\mu^{\nu/2}K_{\nu}\big(\sqrt{\mu\lambda}\big)}\cdot
x^{\nu-1}\cdot\exp\Big\{-\frac12\Big(\frac{\mu}{x}+\lambda
x\Big)\Big\},\ \ \ x>0.%\eqno(33)
$$
Here $\mu>0$, $\lambda\ge0$ if $\nu<0$; $\mu>0$, $\lambda>0$ if
$\nu=0$ and $\mu\ge0$, $\lambda>0$ if $\nu>0$, $K_{\nu}(z)$ is the
modified Bessel function of the third kind with index $\nu$,
$K_{\nu}(z)=\frac12\int_{0}^{\infty}y^{\nu-1}\exp\big\{-\frac{z}{2}\big(y+\frac1y\big)\big\}dy$,
$z\in\mathbb{C}$, $\mathrm{Re}\,z>0$. The corresponding distribution
function will be denoted $P_{GIG}(x;\nu,\mu,\lambda)$.

The class of generalized inverse Gaussian distributions is rather
rich and contains, in particular, both distributions with
exponentially decreasing tails (gamma-distribution ($\mu=0$,
$\nu>0$)), and distributions whose tails demonstrate power-type
behavior (inverse gamma-distribution ($\lambda=0$, $\nu<0$), inverse
Gaussian distribution ($\nu=-\frac12$) and its limit case as
$\lambda\to0$, the L{\'e}vy distribution (the stable distribution
$G_{\frac12,1}(x)$ with the characteristic exponent equal to
$\frac12$ and concentrated on the nonnegative half-line, the
distribution of the time for the standard Wiener process to hit the
unit level)).

In 1977--78 O. Barndorff-Nielsen \cite{BN1977, BN1978} introduced
the class of {\it generalized hyperbolic distributions} as the class
of special normal variance-mean mixtures. For convenience, we will
use a somewhat simpler parametrization. Let $\alpha\in\r$,
$\sigma>0$. If the generalized hyperbolic distribution function with
parameters $\alpha$, $\sigma$, $\nu$, $\mu$, $\lambda$ is denoted
$P_{GH}(x;\alpha,\sigma,\nu,\mu,\lambda)$, then by definition,
$$
P_{GH}(x;\alpha,\sigma,\nu,\mu,\lambda)=\int_{0}^{\infty}\Phi\Big(\frac{x-\alpha
z}{\sigma\sqrt{z}}\Big)\,p_{GIG}(z;\nu,\mu,\lambda)dz,\ \ \
x\in\r.\eqno(28)
$$
Note that in (28) mixing is carried out simultaneously by both
location and scale parameters, but since these parameters are
directly linked in (28), then actually (28) is a one-parameter
mixture.

It is well-known that all generalized hyperbolic distributions are
infinitely divisible \cite{BN1977, BN1978}.

From Theorem 2, Corollary 1 and Theorem B with the account of the
equivalence of relations (11) and (12) we easily obtain the
following result on the convergence of compound Cox processes to
generalized hyperbolic L{\'e}vy processes.

\smallskip

{\sc Corollary 2.} {\it Assume that the random variables
$\{X_{n,j}\}_{j\ge1}$, $n=1,2,...$, $($the jumps of the compound Cox
process $Q_n(t)$, see $(4))$ possess finite variances and satisfy
conditions $(24)$ with some $k_n\in\mathbb{N}$, $a\in\mathbb{R}$ and
$\sigma>0$. Let the the processes $Q_n(t)$ be lead by non-decreasing
positive L{\'e}vy processes $\Lambda_n(t)$ satisfying condition
$(6)$ with some $\delta\in(0,1]$, $\delta_1\ge\frac12$. Also assume
that condition $(26)$ holds. Then the processes $Q_n(t)$ weakly
converge in the Skorokhod space $\mathcal{D}$ to a generalized
hyperbolic L{\'e}vy process $Q(t)$ such that ${\sf
P}\big(Q(1)<x\big)=P_{GH}(x;a,\sigma,\nu,\mu,\lambda)$ if and only
if $ {\sf P}(\Lambda_n(1)<k_nx)\tod P_{GIG}(x;\nu,\mu,\lambda)$ with
the same $k_n$, $\nu$, $\mu$ and $\lambda$. }

\smallskip

The class of generalized gamma (GG) distributions was introduced in
the paper \cite{Stacy1962}. Any representative of this class is
defined by the probability density
$$
p_{GG}(x;\nu,\kappa,\delta)=\frac{|\nu|}{\delta^{k\nu}\Gamma(\kappa)}x^{\kappa\nu-1}\exp\Big\{-\Big(\frac{x}{\delta}\Big)^{\nu}\Big\},\
\ \  x\ge0,
$$
where the parameters $\nu\in\r$, $\kappa,\,\delta\in\r^+$
respectively determine the {\it power, shape and scale} of the
generalized gamma distribution. Here $\Gamma(\kappa)$ is Euler's
gamma-function.

The family of GG-distributions contains practically all most popular
absolutely continuous distributions concentrated on $\r^+$. In
particular, this class contains (i) the gamma-distribution
($\nu=1)$, (ii) the exponential distribution ($\nu=1,\, \kappa=1)$,
(iii) the Erlang distribution ($\nu=1,\, \kappa\in\mathbb{N})$, (iv)
the chi-square distribution ($\nu=1,\, \delta=2)$, (v) the Nakagami
distribution ($\nu=2$), (vi) the half-normal distribution (the
distribution of the absolute value of the standard normal random
variable) ($\nu=2,\,\kappa=\frac12$), (vii) the Rayleigh
distribution ($\nu=2,\,\kappa=1$), (viii) the chi-distribution
($\nu=2,\,\delta=\sqrt{2}$), (ix) the Maxwell distribution
($\nu=2,\,\kappa=3/2$), (x) the Weibull--Gnedenko distribution
($\kappa=1$), (xi) the inverse gamma-distribution ($\nu=-1$), (xii)
the L{\'e}vy distribution, (xiii) the lognormal distribution
($\kappa\to\infty$).

Unlike generalized inverse Gaussian laws, generalized gamma laws
contain the distributions with exponential power type of decrease of
the tail where the posiltive exponent power may be arbitrary.

In \cite{KorolevSokolov2012} the family of generalized variance
gamma (GVG) distributions was introduced as the family of normal
variance-mean mixtures of the form
$$
P_{GVG}(x;a,\sigma,\nu,\kappa,\delta)=\int_{0}^{\infty}\Phi\bigg(\frac{x-au}{\sigma\sqrt{u}}\bigg)
p_{GG}(u;\nu,\kappa,\delta)du,\ \ \ x\in\r,
$$

Both classes, GH and GVG distributions, contain $(a)$ symmetric and
non-symmetric (skew) Student distributions (including Cauchy
distribution), to which in (32) there correspond inverse gamma
mixing distributions; $(b)$ variance gamma (VG) distributions)
(including symmetric and non-symmetric Laplace distributions), to
which in (32) there correspond gamma mixing distributions; $(c)$
normal$\backslash\!\backslash$inverse Gaussian (NIG) distributions
to which in (32) there correspond inverse Gaussian mixing
distributions. However, the class of GVG laws includes normal
mixtures with tails decreasing as $e^{-|x|^{\nu}}$ with $0<\nu<1$,
which are not included in GH-distributions. This type of tail
behavior is of great practical interest.

Not all of GVG distributions are infinitely divisible. But as
concerns the important case of Weibull-type decreasing tails
mentioned above, it turns out that for $\nu\in(0,1]$ these laws are
infinitely divisible. Recall that in the notation introduced above
the Weibull--Gnedenko distribution corresponds to the GG-density
$p_{GG}(x;\nu,1,\delta)$.

\smallskip

{\sc Lemma 5.} {\it If $\nu\le 1$, then the Weibull--Gnedenko
distribution is infinitely divisible.}

\smallskip

{\sc Proof.} As was shown in the paper \cite{KorolevSokolov2014},
the Weibull--Gnedenko distribution with $\nu\le 1$ is mixed
exponential. But in \cite{Goldie1967} it was proved that all mixed
exponential laws are infinitely divisible. The lemma is proved.

\smallskip

According to lemma 5, if $\nu\le 1$, then a L{\'e}vy process $U(t)$,
$t\ge0$, can be defined so that ${\sf
P}\big(U(1)<x\big)=P_{GG}(x;\nu,1,\delta)$, $x\in\mathbb{R}$. Such a
process will be called the {\it L{\'e}vy--Weibull process}. Let
$W(t)$ be a standard Wiener process independent of $U(t)$.

\smallskip

{\sc Corollary 3.} {\it Assume that the random variables
$\{X_{n,j}\}_{j\ge1}$, $n=1,2,...$, $($the jumps of the compound Cox
process $Q_n(t)$, see $(4))$ possess finite variances and satisfy
conditions $(24)$ with some $k_n\in\mathbb{N}$, $a=0$ and
$\sigma>0$. Let the the processes $Q_n(t)$ be lead by non-decreasing
positive L{\'e}vy processes $\Lambda_n(t)$ satisfying condition
$(6)$ with some $\delta\in(0,1]$, $\delta_1\ge\frac12$. Also assume
that condition $(26)$ holds. Then the processes $Q_n(t)$ weakly
converge in the Skorokhod space $\mathcal{D}$ to a subordinated
Wiener process $W(U(t))$ with the subordinator $U(t)$ being the
L{\'e}vy--Weibull process with $\nu\le1$ if and only if $ {\sf
P}(\Lambda_n(1)<k_nx)\tod P_{GG}(x;\nu,1,\delta)$ with the same
$k_n$. }

\renewcommand{\refname}{References}


\begin{thebibliography}{99}

\bibitem{BN1977}
{\it Barndorff-Nielsen O. E.} Exponentially decreasing distributions
for the logarithm of particle size // Proc. Roy. Soc. London.
Ser.~A, 1977. Vol.~A(353). P.~401--419.

\bibitem{BN1978}
{\it Barndorff-Nielsen O. E.} Hyperbolic distributions and
distributions of hyperbolae // Scand. J. Statist., 1978. Vol.~5.
P.~151--157.

\bibitem{BN1979}
{\it Barndorff-Nielsen O. E.} Models for non-Gaussian variation,
with applications to turbulence // Proc. Roy. Soc. London. Ser.~A,
1979. Vol.~A(368). P.~501--520.

\bibitem{BN1982}
{\it Barndorff-Nielsen O. E., Kent J., S{\o}rensen M.} Normal
variance-mean mixtures and z-distributions // International
Statistical Review, 1982. Vol. 50. No. 2. P. 145--159.

\bibitem{BeningKorolev2002} {\it Bening V., Korolev V.} Generalized
Poisson Models and Their Applications in Insurance and Finance. --
Utrecht, VSP, 2002.

\bibitem{Bertoin1996} {\it Bertoin J.} L{\' e}vy Processes. Cambridge Tracts in
Mathematics, Vol. 121. -- Cambridge: Cambridge University Press,
1996.

\bibitem{Billingsley1968} {\it Billingsley P.} Convergence of Probability
Measures. -- New York: Wiley, 1968.

\bibitem{EmbrechtsMaejima2002} {\it Embrechts P., Maejima M.}
Selfsimilar Processes. -- Princeton: Princeton University Press,
2002.

\bibitem{GnedenkoFahim1969} {\it Gnedenko B. V., Fahim H.} On a transfer
theorem // Soviet Math. Dokl., 1969. Vol. 187. No. 1. P. 15--17.

\bibitem{GnedenkoKorolev1996} {\it Gnedenko B. V., Korolev V. Yu.} Random Summation: Limit
Theorems and Applications. -- Boka Raton: CRC Press, 1996.

\bibitem{Goldie1967} {\it Goldie C. M.} A class of infinitely divisible distributions //
Math. Proc. Cambridge Philos. Soc., 1967. Vol. 63. P. 1141-1143.

\bibitem{Grandell1976} {\it Grandell J.} Doubly Stochastic Poisson
Processes. Lecture Notes Mathematics, Vol. 529. --
Berlin--Heidelberg--New York: Springer, 1976.

\bibitem{gz97} {\it Granovsky B. L., Zeifman A. I.} The decay function of
nonhomogeneous birth and death processes, with application to
mean-field models // Stochastic Process. Appl., 1997, Vol. 72. P.
105--120.

\bibitem{Huang2012} {\it Huang H., Kercheval A. N.} A generalized
birth–death stochastic model for high-frequency order book dynamics
// Quantitative Finance, 2012. Vol. 12. No. 4. P. 547--557.

\bibitem{JacodShiryaev2003} {\it Jacod J., Shiryaev A. N.}. Limit theorems for stochastic
processes. 2nd edition. Volume 288 of Grundlehren der Mathematischen
Wissenschaften [Fundamental Principles of Mathematical Sciences]. --
Berlin: Springer-Verlag, Berlin, 2003.

\bibitem{Kashcheev2001} {\it Kashcheev D. E.} Modeling the Dynamics of
Financial Time Series and Evaluation of Derivative Securities. PhD
Thesis. -- Tver: Tver State University, 2001 (in Russian).

\bibitem{Korolev1997} {\it Korolev V. Yu.} On convergence of the
distributions of random sums of independent random variables to
stable laws // Theory Probab. Appl., 1997. Vol. 42. No. 4. P.
818--820.

\bibitem{Korolev2000} {\it Korolev V. Yu.} Asymptotic properties of extrema of
compound Cox processes and their application to some problems of
financial mathematics // Theory Probab. Appl., 2000. Vol. 45. No. 1.
P. 182--194.

\bibitem{Korolev2011} {\it Korolev V. Yu.} Probabilistic and Statistical Methods For Decomposition
of Volatility of Chaotic Processes. -- Moscow: Moscow State
University Publishing House, 2011 (in Russian).

\bibitem{KorolevBeningShorgin2011} {\it Korolev V. Yu., Bening V. E., Shorgin S. Ya.}
Mathematical Foundation of Risk Theory. 2nd ed. -- Moscow:
FIZMATLIT, 2011 (in Russian).

\bibitem{Korolev2013} {\it Korolev V. Yu.} Generalized hyperbolic laws
as limit distributions for random sums // Theory of Probability and
Its Applications, 2013. Vol. 58. No. 1. P. 117-–132.

\bibitem{KorolevSokolov2012} {\it Korolev V. Yu., Sokolov I. A.} Skew Student
distributions, variance gamma distributions and their
generalizations as asymptotic approximations // Informatics and Its
Applications, 2012. Vol. 6. No. 1. P. 2--10.

\bibitem{KorolevZaks2013} {\it Zaks L. M., Korolev V. Yu.} Generalized variance
gamma distributions as limit laws for random sums
// Informatics and Its Applications, 2013. Vol. 7. No. 1. P. 105--115.

\bibitem{KorolevZaksZeifman2013a} {\it Korolev V. Yu., Zaks L. M., Zeifman A. I.}
On convergence of random walks generated by compound Cox processes
to L{\'e}vy processes // Statistics and Probability Letters, 2013.
Vol. 83. No. 10. P. 2432--2438.

\bibitem{Korolev_et_al_2013a} {\it Chertok A., Korolev V., Korchagin A.} Modeling
high-frequency non-homogeneous order flows by compound Cox processes
// Journal of Mathematical Sciences,
2015 (to appear). Available at SSRN, January 14, 2014.
http://ssrn.com/abstract=2378975.

\bibitem{KorolevSokolov2014} {Korolev V. Yu., Sokolov I. A.} On conditions for
the convergence of the distributions of extreme order statistics to
the Weibull distribution // Informatics and its Applications, 2014.
Vol. 8. No. 3. P. 2--10.

\bibitem{KorolevChertokKorchaginZeifman2015} {\it Korolev V. Yu, Chertok A. V.,
Korchagin A. Yu, Zeifman A. I.} Modeling high-frequency order flow
imbalance by functional limit theorems for two-sided risk processes
// Applied Mathematics and Computation, 2015. Vol. 253. P. 224--241

\bibitem{Sato1999} {\it Sato K.} L{\' e}vy Processes and
Infinitely Divisible Distributions. -- Cambridge: Cambridge
University Press, 1999.

\bibitem{Shiryaev1999} {\it Shiryaev A. N.} Essentials of Stochastic Finance: Facts,
Models, Theory. -- Singapore: World Scientific. 1999.

\bibitem{Stacy1962} {\it Stacy E. W.} A generalization of the gamma
distribution // Annals of Mathematical Statistics, 1962. Vol.~33.
P.~1187--1192.

\bibitem{Zolotarev1983} {\it Zolotarev V. M.} One-Dimensional Stable
Distributions. -- Providence, R.I.: American Mathematical Society,
1986.

\end{thebibliography}
\end{document}